\documentclass[reqno,a4paper,11pt]{amsart}
\usepackage{mathrsfs}

\usepackage{amssymb,amsmath,amsfonts,amsthm,breqn}
\usepackage{color}

\newtheorem{theorem}{Theorem}

\theoremstyle{remark}
\newtheorem{rmk}{Remark}

\newcommand{\loc}{\rm{loc}}
\newcommand{\R}{\mathbb{R}}
\renewcommand{\div}{{\rm div}\,}

\subjclass[2010]{35Q31, 93D05.}

\keywords{Euler equations; vorticity; continuity; stability; renormalized solutions}

\title[Strong continuity for Euler]
{Strong continuity for the 2D Euler equations}

\author[G.~Crippa]{Gianluca Crippa}
\address{G.C.: Departement Mathematik und Informatik, Universit\"at Basel, 
Spiegelgasse 1,\break 4051 Basel, Switzerland}
\email{gianluca.crippa@unibas.ch}
\author[E.~Semenova]{Elizaveta Semenova}
\address{E.S.: Department of Epidemiology and Public Health, Swiss Tropical and Public Health Institute, University of Basel, 
Socinstrasse 57, 4051 Basel, Switzerland}
\email{elizaveta.semenova@unibas.ch}
\author[S.~Spirito]{Stefano Spirito}
\address{S.S.: GSSI -- Gran Sasso Science Institute, Viale Francesco Crispi 7, 67100 L'Aquila, Italy}
\email{stefano.spirito@gssi.infn.it}

\begin{document}

\maketitle

\begin{abstract}
We prove two results of strong continuity with respect to the initial datum for bounded solutions to the Euler equations in vorticity form. The first result provides sequential continuity and holds for a general bounded solution. The second result provides uniform continuity and is restricted to H\"older continuous solutions. 
\end{abstract}

\section{Introduction}

Let us consider the Euler equations for an incompressible fluid:
\begin{equation}\label{e:euler}
\begin{cases}
\partial_t u + u \cdot \nabla u = - \nabla p \\
\div u = 0 \,.
\end{cases}
\end{equation}
In two space dimensions, the vorticity $\omega = {\rm curl}\, u$ is a scalar and satisfies the continuity equation
\begin{equation}\label{e:vorticity}
\partial_t \omega + \div(u \, \omega) = 0 \,,
\end{equation}
where the velocity $u$ can be recovered from the vorticity using the Biot-Savart law:
\begin{equation}\label{e:biot}
u(t,x) = \frac{1}{2\pi} \int_{\R^2} \frac{(x-y)^\perp}{|x-y|^2} \, \omega(t,y) \, dy = K \ast \omega \,,
\end{equation}
where $K(x) = x^\perp / (2 \pi |x|^2)$ is the Biot-Savart kernel. We refer to \cite{majda,pulvi} for a comprehensive presentation of the existence and uniqueness theory for the Cauchy problem for the Euler equations. For the purposes of the present paper, it is sufficient to mention that \cite{yudovich} provides existence and uniqueness in the class of bounded vorticities if the initial vorticity $\bar \omega$ belongs to $L^1 \cap L^\infty(\R^2)$, while \cite{dpm2} establish existence of a solution of \eqref{e:euler} for initial vorticities in $L^1 \cap L^p(\R^2)$, where $p>1$ (the uniqueness question is however an open problem).

A fundamental question in fluid dynamics is the continuity of the solution with respect to the initial datum. In the context of bounded vorticities, it is not difficult to prove continuity with respect to the Wasserstein norm (a weak norm arising in the theory of optimal mass transportation). The proof is based on the almost-Lipschitz continuity of the fluid trajectories and the continuity estimate involves a double exponential function of the time (see \cite{loeper,pulvi}).

\medskip

More difficult is to prove continuity estimates with respect to strong norms. To the best of the authors' knowledge, the first result in this context is \cite{wanpulvi}, where stability in $L^1$ of a circular vortex patch was proven. Notice that a circular vortex patch is a stationary solution of the Euler equations. This result was generalized (and the proof simplified) in~\cite{siderisvega}. The only known extension to non-stationary solutions involves elliptic vortex patches \cite{wan,tang}, solutions with a rigid geometry that move with constant angular speed. Nothing is known about general vortex patches.

\medskip

In this note we provide two continuity results with respect to strong norms for non-stationary solutions to the Euler equations, without any geometric requirement on the shape of the solutions. Our first theorem holds for general bounded solutions~$\omega$ of the Euler equations and provides strong convergence in space at time $t>0$, provided that the initial data converge strongly.

\begin{theorem}\label{t:punct}
Let $\bar \omega \in L^1 \cap L^\infty(\R^2)$ and let $\{\bar \omega_n\}_n \subset L^1 \cap L^2(\R^2)$ be a sequence with $\bar \omega_n \to \bar \omega$ strongly in $L^2(\R^2)$. Fix $T<\infty$ and for every $n$ let $\omega_n \in C([0,T];L^1 \cap L^2(\R^2))$ be a solution of the Euler equations \eqref{e:vorticity}--\eqref{e:biot} with initial datum $\bar \omega_n$. Then 
$$
\omega_n(t,\cdot) \to \omega(t,\cdot) \qquad \text{strongly in $L^2(\R^2)$, uniformly for $t \in [0,T]$,}
$$
where $\omega$ is the unique solution in $C([0,T];L^1 \cap L^\infty(\R^2))$ of the Euler equations~\eqref{e:vorticity}--\eqref{e:biot} with initial datum $\bar \omega$.
\end{theorem}

First, we want to point out that the continuity in time with values in $L^1(\R^2)\cap L^{2}(\R^2)$ follows from the fact that for any $n$ the the vorticity $\omega_n$ is a renormalized solutions of~\eqref{e:vorticity} and the velocity $u_n$ is in the setting of~\cite{diplio}, see Step~1 in the proof of Theorem~\ref{t:punct}. Moreover, notice that in the above theorem we do not require boundedness of $\bar \omega_n$, therefore the Euler equations with such initial data may have more than one solution. Nevertheless, thanks to the uniqueness for the (limit) problem with bounded initial datum $\bar \omega$, the result holds for any sequence of solutions, and does not require the passage to a subsequence.

The proof relies on the DiPerna-Lions theory of continuity equations with Sobolev velocity field \cite{diplio}, see also \cite{hw} for a general account of this research area. This proof will be presented in \S\ref{s:A}.

\begin{rmk}
If the limit vorticity $\bar \omega$ only belongs to $L^1 \cap L^2(\R^2)$ we only have convergence of a {\em subsequence} of $\omega_n(t,\cdot)$ to {\em some} solution of the limit problem, since the latter has no uniqueness: Step~4 in the proof of Theorem~\ref{t:punct} does not apply. See also~\cite{anna} for a proof of the existence of solutions to the Euler equations via Lagrangian techniques.
\end{rmk}

\begin{rmk} 
If we consider in Theorem~\ref{t:punct} a sequence of initial data $\{\bar\omega_n\}_n \subset L^1 \cap L^p(\R^2)$ converging to $\bar \omega$ strongly in $L^p(\R^2)$, with $p>2$, then it is possible to prove the strong convergence in $L^p(\R^2)$ of $\omega_n(t,\cdot)$ to $\omega(t,\cdot)$. On the other hand, if we relax the integrability assumption on the sequence $\bar \omega_n$ to $L^1 \cap L^p(\R^2)$ for some  $p<2$, our proof breaks down in its full generality. Indeed, a vorticity in $L^p$ advected by a velocity in $W^{1,p}_\loc$ with $p<2$ does not fall in the context of \cite{diplio}, therefore existence of a flow as in Step~1 of the proof of Theorem~\ref{t:punct} is not guaranteed (in fact, if $p<4/3$ equation~\eqref{e:vorticity} does not even make distributional sense). One should consider a sequence of approximate solutions $\omega_n$ that are a priori required to be Lagrangian. Note that \cite{stefano} guarantees that solutions obtained via vanishing viscosity approximation are indeed Lagrangian.
\end{rmk}

The above theorem provides {\em sequential} continuity (with respect to the $L^2$ norm) of the map that associates to the initial datum the solution at time $t$. This continuity property holds at {\em every} bounded solution. The rate of continuity may however depend on the solution itself. Our second result provides uniform continuity with an explicit convergence rate, provided we restrict our attention to slightly more regular solutions to the Euler equations~\eqref{e:vorticity}--\eqref{e:biot}: we require that the (compactly supported) initial data (and therefore the solution at any time) belong to some H\"older class $C^\alpha_c(\R^2)$, where $0<\alpha<1$ is arbitrary.

\begin{theorem} \label{t:unif}
Assume that $\bar \omega_1$, $\bar \omega_2 \in C^\alpha_c(\R^2)$, where $0<\alpha<1$, satisfy $\int_{\R^2} \bar \omega_1 \, dx = \int_{\R^2} \bar \omega_2 \, dx $, and fix $T<\infty$. Let $\omega_1$ and $\omega_2$ be the unique solutions in $C([0,T];L^1 \cap L^\infty(\R^2))$ of the Euler equations \eqref{e:vorticity}--\eqref{e:biot} with initial data $\bar \omega_1$ and $\bar \omega_2$, respectively. Then
\begin{equation}\label{e:thm}
\| \omega_1(t,\cdot) - \omega_2(t,\cdot) \|_{L^2} \leq C e^{ct} \| \bar \omega_1 - \bar \omega_2 \|^\gamma_{L^2}
\qquad \text{for $t\in [0,T]$,}
\end{equation}
where $C$, $c$, and $\gamma$ only depend on $\alpha$, $T$, the norms of $\bar \omega_1$ and $\bar \omega_2$, and the diameter of the supports of $\bar \omega_1$ and $\bar \omega_2$.
\end{theorem}

The proof of this theorem involves an interpolation argument in homogeneous fractional Sobolev spaces. Essentially, the H\"older regularity of the solution allows to ``upgrade'' weak estimates (as in~\cite{loeper,pulvi}) to strong estimates. The proof will be presented in \S\ref{s:B}.

\begin{rmk} 
Inequality \eqref{e:thm} can be extended to $L^p$ norms with $1 \leq p < \infty$, although with a different value for the constants and the exponent. 
\end{rmk}

\subsection*{Acknowledgment} This research has been partially supported by the SNSF grants 140232 and 156112.

\section{Proof of Theorem~\ref{t:punct}}
\label{s:A}

\subsection*{Step 1.}
Let us consider the velocity $u_n$ associated to the vorticity $\omega_n$ as in \eqref{e:biot}. Decomposing the Biot-Savart kernel as $K = K_1 + K_2 = K {\bf 1}_{|x|\leq 1} + K {\bf 1}_{|x| > 1}$ and noting that $K_1 \in L^1(\R^2)$ and $K_2 \in L^\infty(\R^2)$, we obtain with Young's inequality that 
$$
u_n \in L^\infty ([0,T]; L^1(\R^2) + L^\infty (\R^2)) \,.
$$
In particular, formula \eqref{e:biot} is well-defined in this summability context. Moreover, $u_n$ is divergence-free and (by elliptic regularity, since $\omega \in L^\infty([0,T];L^2(\R^2))$) belongs to $L^\infty([0,T] ; W^{1,2}_\loc(\R^2))$.

The bounds above imply that we are in the setting of \cite{diplio}: there exists a unique forward-backward regular Lagrangian flow (i.e., in this context, an incompressible flow defined almost everywhere in space) $X_n = X_n(s,t,x)$ associated to the velocity field $u_n$, and the vorticity $\omega_n$ is transported by such a flow, in the sense that
\begin{equation}\label{e:lag}
\omega_n(t,x) = \bar \omega_n (X_n(0,t,x)) \,.
\end{equation}

\subsection*{Step 2.} From the representation \eqref{e:lag}, together with the convergence of $\bar \omega_n$ to $\bar \omega$, it follows that $\omega_n \in L^\infty([0,T];L^1\cap L^2(\R^2))$ uniformly in $n$. Therefore, along a subsequence we have $\omega_{n(k)} \rightharpoonup w$ weakly* in $L^\infty([0,T];L^1\cap L^2(\R^2))$. Moreover, all the bounds on $u_n$ listed in Step~1 are uniform in $n$. Arguing as in \cite[Theorem~1.2]{dpm2} we find a further subsequence (that we do not relabel) $u_{n(k)}$ converging strongly in $L^2_\loc([0,T] \times \R^2)$ to a limit velocity $v$ (notice that the convergence is strong also with respect to the time: this makes use of Aubin's lemma). One can readily check that $v$ enjoys the same bounds as in Step~1 for the sequence $u_n$ and that the couple $(v,w)$ solves~\eqref{e:vorticity}--\eqref{e:biot}.

\subsection*{Step 3.} We can therefore apply the stability theorem from \cite{diplio} (see also \cite{crelle,jhde} for a purely Lagrangian proof of such a stability theorem, and \cite{abc1,abc2,abc3}, specific to the two-dimensional context). We obtain that the flows $X_{n(k)}$ from Step~1 converge locally in measure in $\R^2$, uniformly in $t$, $s \in [0,T]$, to the unique forward-backward regular Lagrangian flow $X$ associated to the velocity field $v$.

Therefore 
\begin{equation}\label{e:lastlag}
\omega_{n(k)}(t,x) = \bar \omega_{n(k)} (X_{n(k)}(0,t,x)) \to \bar \omega (X(0,t,x)) 
\end{equation}
strongly in $L^2(\R^2)$, uniformly for $t \in [0,T]$ (here one can argue using Lusin's theorem and exploiting the incompressibility of the flows, see for instance the argument in \cite[Propositions~7.2 and~7.3]{jhde}). Hence, the weak limit $w$ of $\omega_n$ defined in Step~2 is in fact a strong limit and coincides with $\bar \omega (X(0,t,x))$.

\subsection*{Step 4.}
The representation in \eqref{e:lastlag} entails that $w$ is a bounded function that solves~\eqref{e:vorticity}--\eqref{e:biot}, therefore by uniqueness it coincides with the solution $\omega$ in the statement of the theorem. By uniqueness of the limit the whole sequence $\omega_n(t,\cdot)$ (and not only the subsequence $\omega_{n(k)}(t,\cdot)$, as in \eqref{e:lastlag}) converges to $\omega(t,\cdot)$. This concludes the proof of Theorem~\ref{t:punct}. \qed

\section{Proof of Theorem~\ref{t:unif}}
\label{s:B}

First of all, we observe that, since $\omega_1$, $\omega_2 \in L^\infty ([0,T]; L^1\cap L^\infty (\R^2))$, the velocities $u_1$ and $u_2$ are uniformly bounded. This in turn implies that $\omega_1$ and $\omega_2$ are compactly supported in space, uniformly for $t \in [0,T]$. In the course of the proof, we will use the notation $\dot{H}^s(\R^2)$ to denote the homogenous Sobolev space in $\R^2$ of real order $s$.

\subsection*{Step 1.} Let us fix $0<\beta<\alpha$. The classical interpolation inequality for homogeneous Sobolev spaces gives
\begin{equation}\label{e:proof1}
\| \omega_1(t,\cdot) - \omega_2 (t,\cdot) \|_{L^2} 
\leq \| \omega_1(t,\cdot) - \omega_2 (t,\cdot) \|^{\frac{\beta}{1+\beta}}_{\dot{H}^{-1}}  
\| \omega_1(t,\cdot) - \omega_2 (t,\cdot) \|^{\frac{1}{1+\beta}}_{\dot{H}^\beta} \,.
\end{equation}
It is known (see for instance \cite[Page 326]{majda}) that H\"older regularity of the vorticity is propagated in time by~\eqref{e:vorticity}, and it is immediate to check that $C^\alpha_c (\R^2) \hookrightarrow H^\beta (\R^2)$ for $\alpha>\beta$. Since moreover $\int_{\R^2} \omega_1(t,\cdot) - \omega_2(t,\cdot) \, dx = 0$ for all times, the second factor in the right hand side of \eqref{e:proof1} is bounded by a constant uniformly in time.

\subsection*{Step 2.} Since $\omega_i \in C^\alpha_c (\R^2)$ uniformly in time, by elliptic regularity $u_i \in {\rm Lip}_\loc (\R^2)$ uniformly in time (see for instance \cite[Page 327]{majda}). A standard $L^2$-energy estimate for the difference of the equations \eqref{e:euler} for $u_1$ and $u_2$ then implies
\begin{equation}\label{e:proof2}
\| u_1(t,\cdot) - u_2(t,\cdot) \|_{L^2} \leq e^{ct} \| \bar u_1 - \bar u_2 \|_{L^2}\,.
\end{equation}
Notice that $u_1-u_2$ is globally in $L^2(\R^2)$ since $\omega_1-\omega_2$ has zero integral (see~\cite[Page 321]{majda}).

\subsection*{Step 3.} Since the velocity is divergence-free one can check (for instance, passing in Fourier variables) that
\begin{equation}\label{e:proof3}
\| \omega_1 (t,\cdot) - \omega_2 (t,\cdot) \|_{\dot{H}^{-1}} = \| u_1(t,\cdot) - u_2(t,\cdot) \|_{L^2} \,.
\end{equation}
Moreover, since $\bar \omega_1- \bar \omega_2$ has zero integral, elliptic regularity implies the global estimate
\begin{equation}\label{e:proof4}
\| \bar u_1 - \bar u_2 \|_{L^2} \leq C \| \bar \omega_1 - \bar \omega_2 \|_{L^2}\,.
\end{equation}

\subsection*{Step 4.} Combining \eqref{e:proof1}--\eqref{e:proof4} we obtain \eqref{e:thm}. \qed

\end{document}